\documentclass[11pt,a4paper]{article}

\usepackage{inputenc}
\usepackage{amsmath}
\usepackage{bm}
\usepackage{bbold}
\usepackage{amsthm}

\usepackage{hyperref}

\title{Recursive Equations Based Models\\ of Queueing Systems\thanks{Proc. 1994 SCS Europ. Simulation Symp., Istanbul, Turkey, Oct. 9-12, 1994 / Ed. by A. R. Kaylan, A. Lehmann, T. I. Oren. P. 252--256.}
}

\author{Nikolai K. Krivulin\thanks{Faculty of Mathematics and Mechanics, St.~Petersburg State University, Bibliotechnaya Sq.2, Petrodvorets, St.~Petersburg, 198904 Russia}}

\date{}

\begin{document}

\maketitle

\begin{abstract}
An overview of the recursive equations based models and their applications in
simulation based analysis and optimization of queueing systems is given. These
models provide a variety of systems with a convenient and unified
representation in terms of recursions for arrival and departure times of
customers, which involves only the operations of maximum, minimum, and
addition.
\end{abstract}

\section{Introduction}

As a representation of dynamics of queueing systems, recursive equations have
been introduced by Lindley 1952 in his investigation of the $  G/G/1  $
queue. The representation has proved to be useful in both analytical study and
simulation of queues, and was extended to cover a variety of queueing systems
including open and closed tandem single-server queues with both infinite and
finite buffers, the $  G/G/m  $ system, and queueing networks.

The recursive equations were originally expressed in terms of the waiting
times of customers (\cite{Lindley52,Kiefer55}). Equations
following this classical approach remain traditional in the queueing theory,
one can find them in many of the recent works devoted mainly to theoretical
aspects of the investigation of queueing systems (see, e.g., \cite{Kalashnikov94}).

In the last few years, another representation based on recursions for the
arrival and departure times of customers has gained acceptance in works
dealing with the simulation study of queueing systems and its related fields
including performance evaluation and sensitivity analysis. The items of our
list of references, other than those cited above, can serve as an
illustration. Although these equations may be readily derived from those of
the classical type, they offer a more convenient and unified way of
representing dynamics of queueing systems as well as their performance
measures.

The purpose of this paper is to give a brief overview of the recursive
equations and their applications in simulation based analysis and optimization
of queueing systems. The subsequent sections present the equations expressed
in terms of the arrival and departure times, which describe the systems most
commonly encountered in studies of queues, and discuss the representation of
performance measures, associated with these queueing system models.
Applications of the models to the development of simulation algorithms as well
as to the analysis of system performance measures and estimation of their
sensitivity are also outlined. Finally, limitations on the use of the models
are briefly discussed.

\section{Models of Queueing Systems}

Most of the models appearing in this section actually present single-server
systems which can have both finite and infinite buffers, and operate according
to the first-come, first-served (FCFS) queueing discipline. Also included are
the equations representing the $  G/G/m  $ system, and a rather general
model of a queueing network with a deterministic routing mechanism.

\subsection{The $  G/G/1  $ Queue}

We start with this model which provides the basis for representing more
complicated queueing systems. The $  G/G/1  $ system consists of a server
and a buffer with infinite capacity (Fig.1). Once a customer arrives into the
system, he occupies the server provided that it is free. If the customer finds
the server busy, he is placed into the buffer, and starts waiting to be
served. The queue discipline in the system is presumed to be FCFS.
\begin{figure}[hhh]
\setlength{\unitlength}{1.5mm}
\begin{center}
\begin{picture}(25,10)

\put(5,6){$A_{k}$}
\put(19,6){$D_{k}$}

\put(0,4){\line(1,0){4}}
\put(0,0){\line(1,0){4}}
\put(4,4){\line(0,-1){4}}

\put(10,4){\line(1,0){4}}
\put(10,0){\line(1,0){4}}
\put(14,4){\line(0,-1){4}}

\put(4,2){\vector(1,0){6}}
\put(18,2){\vector(1,0){6}}

\put(16,2){\circle{3}}

\end{picture}
\end{center}
\caption{The $  G/G/1  $ queue.}
\end{figure}
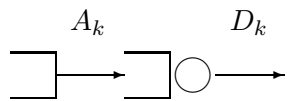

For the $  G/G/1  $ queue, we denote the interarrival time between the $k$th
customer and his predecessor by $  \alpha_{k} $, and the service time of the
$k$th customer by $  \tau_{k} $. Furthermore, let $  A_{k}  $ be the $k$th
arrival epoch to the queue, and $  D_{k}  $ be the $k$th departure epoch
from the queue, $  k=1,2, \ldots $. Provided that the system starts operating
at time zero, it is convenient to set $  A_{0} \equiv 0  $ and
$  D_{0} \equiv 0 $. One may now describe the dynamics of the $  G/G/1  $
queue as
\begin{eqnarray*}
A_{k} & = & A_{k-1} + \alpha_{k} \\
D_{k} & = & (A_{k} \vee D_{k-1}) + \tau_{k},
\end{eqnarray*}
where $  \vee  $ stands for the maximum operator, $  k=1,2, \ldots $.

\subsection{Tandem Systems of Single-Server Queues}

Consider a system of $  N  $ single-server queues with infinite buffers,
operating in tandem as shown in Fig.2.
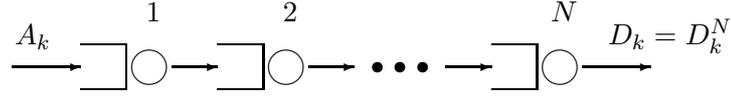
\begin{figure}[hhh]
\setlength{\unitlength}{1.5mm}
\begin{center}
\begin{picture}(60,10)

\put(1,4){$A_{k}$}
\put(53,4){$D_{k} = D_{k}^{N}$}

\put(12.5,6){$1$}
\put(24.5,6){$2$}
\put(48,6){$N$}

\put(7,4){\line(1,0){4}}
\put(19,4){\line(1,0){4}}
\put(43,4){\line(1,0){4}}

\put(7,0){\line(1,0){4}}
\put(19,0){\line(1,0){4}}
\put(43,0){\line(1,0){4}}

\put(11,4){\line(0,-1){4}}
\put(23,4){\line(0,-1){4}}
\put(47,4){\line(0,-1){4}}

\put(1,2){\vector(1,0){6}}
\put(15,2){\vector(1,0){4}}
\put(27,2){\vector(1,0){4}}
\put(39,2){\vector(1,0){4}}
\put(51,2){\vector(1,0){6}}

\put(13,2){\circle{3}}
\put(25,2){\circle{3}}
\put(49,2){\circle{3}}

\multiput(33,2)(2,0){3}{\circle*{1}}

\end{picture}
\end{center}
\caption{Single-server queues in tandem.}
\end{figure}

Each customer that arrives into this system has to pass through all the queues
so as to occupy consecutively every servers from $  1  $ to $  N $, and
then leave the system. We suppose that upon his service completion at a queue,
the customer arrives into the next queue immediately.

To set up the recursive equations representing the system in a convenient way,
let us introduce the symbols $  D_{k}^{n}  $ and $  \tau_{k}^{n}  $
respectively for the departure and service times of the $k$th customer at
queue $  n $. However, we maintain the symbols $  A_{k}  $ and
$  D_{k} = D_{k}^{N}  $ to denote the $k$th arrival and departure epochs for
the whole system. With these notations, the equations are written as
(Shanthikumar and Yao 1989a; and Chen and Chen 1990)
\begin{eqnarray*}
D_{k}^{1} & = & (A_{k} \vee D_{k-1}^{1}) + \tau_{k}^{1} \\
D_{k}^{n} & = & (D_{k}^{n-1} \vee D_{k-1}^{n}) + \tau_{k}^{n},
\quad
n = 2, \ldots, N.
\end{eqnarray*}

\subsubsection{Closed Tandem Systems}

Suppose that in the above tandem system all
the customers after their service completion at the $N$th server return to the
$1$st queue for the next cycle of service (see Fig.3).
\begin{figure}[hhh]
\setlength{\unitlength}{1.5mm}
\begin{center}
\begin{picture}(56,18)

\put(11.5,12){$1$}
\put(23.5,12){$2$}
\put(47,12){$N$}

\put(6,10){\line(1,0){4}}
\put(18,10){\line(1,0){4}}
\put(42,10){\line(1,0){4}}

\put(6,6){\line(1,0){4}}
\put(18,6){\line(1,0){4}}
\put(42,6){\line(1,0){4}}

\put(10,10){\line(0,-1){4}}
\put(22,10){\line(0,-1){4}}
\put(46,10){\line(0,-1){4}}

\put(0,8){\vector(1,0){6}}
\put(14,8){\vector(1,0){4}}
\put(26,8){\vector(1,0){4}}
\put(38,8){\vector(1,0){4}}
\put(50,8){\vector(1,0){6}}

\put(0,8){\line(0,-1){8}}
\put(56,8){\line(0,-1){8}}
\put(0,0){\line(1,0){56}}

\put(12,8){\circle{3}}
\put(24,8){\circle{3}}
\put(48,8){\circle{3}}

\multiput(32,8)(2,0){3}{\circle*{1}}

\put(6,2){$K_{1}$}
\put(18,2){$K_{2}$}
\put(42,2){$K_{N}$}

\end{picture}
\end{center}
\caption{A closed tandem system of single-server queues.} 
\end{figure}
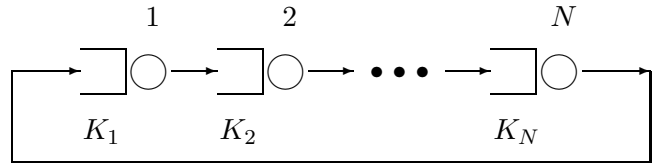

Furthermore, we assume that at the initial time, there are $  K_{n} $,
$  0 \leq K_{n} < \infty $, customers in the buffer of server $  n $.
Assuming $  D_{k}^{n} = -\infty  $ for all $  k < 0  $ and
$  n=1, \ldots N $, one can represent the closed system in the form (see,
e.g., \cite{Greenberg91,Krivulin94b})
\begin{eqnarray*}
D_{k}^{1} & = & (D_{k-K_{1}}^{N} \vee D_{k-1}^{1}) + \tau_{k}^{1}  \\
D_{k}^{n} & = & (D_{k-K_{n}}^{n-1} \vee D_{k-1}^{n}) + \tau_{k}^{n},
\quad
n = 2, \ldots, N.
\end{eqnarray*}

\subsubsection{Tandem Queues with Finite Buffers}

We now turn to the discussion of
the system of queues which provide only a limited number of places in their
buffers for customers waiting for service. In such a system, if the buffer at
a server has finite capacity, the preceding server may be blocked according to
one of the blocking rules. In this paper we shall restrict ourselves to
{\em manufacturing} blocking and {\em communication} blocking which are more
frequent in practice.

Consider a system of $  N  $ queues, depicted in Fig.4. We denote the
capacity of the buffer at server $  n  $ by $  B_{n} $,
$  0 \leq B_{n} \leq \infty $, $  n=2,3,\ldots,N $. As the input buffer of
the system, the buffer of the $1$st server is assumed to be infinite.
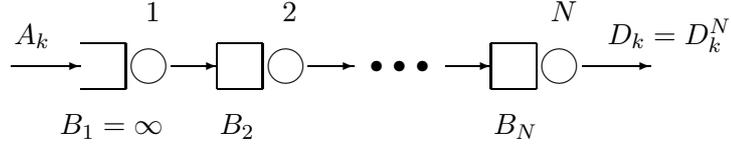
\begin{figure}[hhh]
\setlength{\unitlength}{1.5mm}
\begin{center}
\begin{picture}(56,15)

\put(0,8){$A_{k}$}
\put(52,8){$D_{k} = D_{k}^{N}$}

\put(11.5,10){$1$}
\put(23.5,10){$2$}
\put(47,10){$N$}

\put(6,8){\line(1,0){4}}
\put(18,8){\line(1,0){4}}
\put(42,8){\line(1,0){4}}

\put(6,4){\line(1,0){4}}
\put(18,4){\line(1,0){4}}
\put(42,4){\line(1,0){4}}

\put(18,8){\line(0,-1){4}}
\put(42,8){\line(0,-1){4}}

\put(10,8){\line(0,-1){4}}
\put(22,8){\line(0,-1){4}}
\put(46,8){\line(0,-1){4}}

\put(0,6){\vector(1,0){6}}
\put(14,6){\vector(1,0){4}}
\put(26,6){\vector(1,0){4}}
\put(38,6){\vector(1,0){4}}
\put(50,6){\vector(1,0){6}}

\put(12,6){\circle{3}}
\put(24,6){\circle{3}}
\put(48,6){\circle{3}}

\multiput(32,6)(2,0){3}{\circle*{1}}

\put(4,0){$B_{1} = \infty$}
\put(18,0){$B_{2}$}
\put(42,0){$B_{N}$}

\end{picture}
\end{center}
\caption{Tandem single-server queues with finite buffers.} 
\end{figure}

Let us first suppose that the system operates according to the manufacturing
blocking rule. Under this type of blocking, if a customer upon completion of
his service at server $  n  $ sees the buffer of the $(n+1)$st server full,
he cannot unoccupy the $n$th server until the next server provides a free
space in its buffer. Since buffers become free as customers are called forward
for service, the $n$th server is unoccupied as soon as the $(n+1)$st server
completes its current service to initiate the service of the next customer. It
is not difficult to understand that the departure of the $k$th customer from
server $  n  $ occurs not earlier than that of the $(k-B_{n+1}-1)$st
customer from server $  n+1 $. Taking into account this condition, one may
represent the equations as (\cite{Shanthikumar89a,Chen90})
\begin{eqnarray*}
D_{k}^{1} & = & ((A_{k} \vee D_{k-1}^{1}) + \tau_{k}^{1})
                                          \vee D_{k-B_{2}-1}^{2} \\
D_{k}^{n} & = & ((D_{k}^{n-1} \vee D_{k-1}^{n}) + \tau_{k}^{n})
                                          \vee D_{k-B_{n+1}-1}^{n+1},
\quad
n = 2, \ldots, N-1  \\
D_{k}^{N} & = & (D_{k}^{N-1} \vee D_{k-1}^{N}) + \tau_{k}^{N}.
\end{eqnarray*}

The communication blocking rule requires from a server not to initiate the
service of a customer if the buffer of the next server is full. In this case,
the server remains unavailable until the current service at the next server
is completed. For the system with communication blocking, we have (Chen and 
Chen 1990)
\begin{eqnarray*}
D_{k}^{1} & = & (A_{k} \vee D_{k-1}^{1} \vee D_{k-B_{2}-1}^{2})
                                           + \tau_{k}^{1} \\
D_{k}^{n} & = & (D_{k}^{n-1} \vee D_{k-1}^{n} \vee D_{k-B_{n+1}-1}^{n+1})
                                           + \tau_{k}^{n},
\quad
n = 2, \ldots, N-1  \\
D_{k}^{N} & = & (D_{k}^{N-1} \vee D_{k-1}^{N}) + \tau_{k}^{N}.
\end{eqnarray*}

\subsection{$G/G/m  $ Queues}

Equations representing the $  G/G/m  $ queue (Fig.5) as recursions for the
waiting times of customers have been first introduced by Kiefer and Wolfowitz
1955 \cite{Kiefer55}. These recursive equations were expressed in general terms rather than in
an explicit form similar to those presented above.
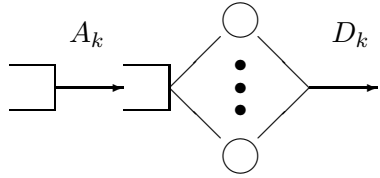
\begin{figure}[hhh]
\setlength{\unitlength}{1.5mm}
\begin{center}
\begin{picture}(35,15)

\put(5,10){$A_{k}$}
\put(28,10){$D_{k}$}

\put(0,8){\line(1,0){4}}
\put(0,4){\line(1,0){4}}
\put(4,8){\line(0,-1){4}}

\put(10,8){\line(1,0){4}}
\put(10,4){\line(1,0){4}}
\put(14,8){\line(0,-1){4}}

\put(14,6){\line(1,1){4.5}}
\put(14,6){\line(1,-1){4.5}}

\put(26.3,6){\line(-1,1){4.5}}
\put(26.3,6){\line(-1,-1){4.5}}

\put(4,6){\vector(1,0){6}}
\put(26.3,6){\vector(1,0){6}}

\put(20.25,12){\circle{3}}
\put(20.25,0){\circle{3}}
\multiput(20.25,4)(0,2){3}{\circle*{1}}

\end{picture}
\end{center}
\caption{The $  G/G/m  $ queue.} 
\end{figure}

To represent the equations for the $  G/G/m  $ queue,
$  1 \leq m < \infty $, in terms of the arrival and departure times of
customers, let us further insert the symbol $  C_{k}  $ for the service
completion time of the customer which is the $k$th to arrive into the system.
Note that in multi-server queues the $k$th departure time and the completion
time of the $k$th customer may not coincide as contrasted to the
$  G/G/1  $ queue which does not recognize them.

Now we may describe the dynamics of the system through the equations proposed
in \cite{Krivulin94a} (a similar representation in terms of waiting times can be
found in \cite{Tsitsiashvili76})
\begin{eqnarray*}
A_{k} & = & A_{k-1} + \alpha_{k} \\
C_{k} & = & (A_{k} \vee D_{k-m}) + \tau_{k} \\
D_{k} & = & \bigwedge_{ 1 \leq j_{1} < \ldots < j_{k} \leq k+m-2 }
( C_{j_{1}} \vee \cdots \vee C_{j_{k}} ) \wedge C_{k+m-1},
\end{eqnarray*}
where $  \wedge  $ signifies the minimum operator. Note that with
$  m=1  $ the above set of equations is reduced to that of the $  G/G/1  $
queue.

\subsection{Networks with Deterministic Routing}

We complete this section by presenting a rather general model of a closed
queueing network with deterministic routing described in \cite{Krivulin90,Krivulin93} (see also a similar model in \cite{Shanthikumar89b}). Let us
first consider a network consisting of $  N  $ single-server nodes. In each
node there are a server and an infinite buffer in which customers are placed
at their arrival so as to wait for service if it cannot be initiated
immediately. After his service completion at one node, each customer goes to
another node chosen according to the routing procedure defined as follows.
For the network, we assume that a matrix
$$
{\bf S} = \left( \begin{array}{ccccc}
                     s_{11} & s_{12} & \cdots & s_{1k} & \cdots \\
                     s_{21} & s_{22} & \cdots & s_{2k} & \cdots \\
	             \vdots & \vdots &        & \vdots &        \\
                     s_{N1} & s_{N2} & \cdots & s_{Nk} & \cdots
                 \end{array}
          \right)
$$
is given, $  s_{nk}  $ determines the next node to be visited by the
customer who is the $k$th to depart from node $  n $,
$  s_{nk} \in \{1, \ldots, N\} $; $  n=1, \ldots, N $; $  k=1,2,\ldots $.
It is also assumed that at the initial time, all servers are free, and there
are $  K_{n} $, $  0 \leq K_{n} \leq \infty $, customers in the buffer at
node $  n $.

For node $  n $, we denote the $k$th arrival and departure epochs
respectively by $  A_{k}^{n}  $ and $  D_{k}^{n} $, and the service time of
the customer who is the $k$th to arrive by $  \tau_{k}^{n} $. Furthermore,
let us introduce the set
$$ \bm{D}_{n} =
       \{D_{i}^{k} |   s_{ik} = n;   i=1,\ldots,N;   k=1,2,\ldots \}, $$
which is constituted by the departure times of the customers who have to go to
node $  n $. Finally, we may represent the network by means of the equations
\begin{eqnarray*}
D_{k}^{n} & = & (A_{k}^{n} \vee D_{k-1}^{n}) + \tau_{k}^{n} \\
A_{k}^{n} & = & \left\{ \begin{array}{ll}
			  0,                     & \mbox{if $k \leq K_{n}$} \\
			  {\cal A}_{k-K_{n}}^{n}, & \mbox{otherwise,}
			\end{array}
		\right.
\end{eqnarray*}
where $  {\cal A}_{k}^{n}  $ is the arrival time of the customer which is
the $k$th to arrive into node $  n  $ after his service at any node of the
network. In other words, the symbol $  {\cal A}_{k}^{n}  $ differs from
$  A_{k}^{n}  $ in that it relates only to the customers really arriving
into node $  n $, and does not to those occurring in this node at the initial
time. It is defined as
$$ {\cal A}_{k}^{n} =
\bigwedge_{ \{ D_{1},\ldots, D_{k} \} \subset \bm{D}_{n} }
( D_{1} \vee \cdots \vee D_{k} ),
$$
where minimum is taken over all $k$-subsets of the set $  \bm{D}_{n} $.

It is easy to understand how tandem queueing systems with infinite buffers
may be represented as networks like that just described. Moreover, changing
the first one from the above equations, one can readily extend the model to
cover networks with nodes which may have many servers. These servers may
operate both in tandem and in parallel, and even form a network themselves.

\section{Performance Measures}

We are now in a position to show how the performance measures which one
normally chooses in the analysis of queueing systems may be represented on the
basis of the models described above. We start with presenting sample
performance measures associated with the systems under consideration, and then
briefly discuss the evaluation of system performance measures in the general
case.

\subsection{Networks with Single-Server Nodes}

Suppose that we observe the network until the $K$th service completion at node
$  n $, $  K=1,2, \ldots $; $  1 \leq n \leq N $. As sample performance
measures for node $  n  $ in the observation period, one can consider the
following average quantities (\cite{Chen90,Krivulin90,Krivulin93}):
\begin{eqnarray*}
S_{K}^{n} &=& \sum_{k=1}^{K} (D_{k}^{n}-A_{k}^{n})/K,
\quad
\text{the total time of one customer};
\\
W_{K}^{n} &=& \sum_{k=1}^{K} (D_{k}^{n}-A_{k}^{n}-\tau_{k}^{n})/K; 
\quad
\text{the waiting time of one customer};
\\
T_{K}^{n} &=& K/D_{K}^{n},
\quad
\text{the throughput rate of the node};
\\
U_{K}^{n} &=& \sum_{k=1}^{K} \tau_{k}^{n}/D_{K}^{n},
\quad
\text{the utilization of the server};
\\
J_{K}^{n} &=& \sum_{k=1}^{K} (D_{k}^{n}-A_{k}^{n})/D_{K}^{n},
\quad
\text{the number of customers at the node};
\\
Q_{K}^{n} &=& \sum_{k=1}^{K} (D_{k}^{n}-A_{k}^{n}-\tau_{k}^{n})/D_{K}^{n},
\quad
\text{the queue length at the node}.
\end{eqnarray*}

Note that, assuming the service times $  \tau_{k}^{n}  $ to be given, one
can express these measures in closed form only in terms of these times,
involving arithmetic operations, and the operations of maximum and minimum
(\cite{Krivulin93}).

\subsection{Tandem Systems of Queues}

Since tandem systems can be considered as networks with deterministic routing,
the above sample performance measures are also suited to the tandem systems.
In addition to these measures which are actually server related performance
criteria, for a tandem system with $  N  $ servers we may define customer
related performance measures \cite{Chen90}
\begin{eqnarray*}
S_{K} &=& \sum_{k=1}^{K} (D_{k} - A_{k})/K,
\quad
\text{the average system time of one customer};
\\
W_{K} &=& \sum_{k=1}^{K} \left( D_{k} - A_{k} - \sum_{n=1}^{N} \tau_{k}^{n} \right)/K,
\quad
\text{the average waiting time of one customer}.
\end{eqnarray*}

Finally, there are sample performance measures inherent only in the systems
with finite buffers. As an example, the average idle time of a server, say
server $  n $, can be considered. This measure is written in the same form
for both the manufacturing and communication blocking rules as
$$
I_{K}^{n} = \sum_{k=1}^{K} \left(D_{k}^{n}-(D_{k}^{n-1}\vee D_{k-1}^{n})
            - \tau_{k}^{n} \right)/K.
$$

\subsection{Multi-Server Queues}

Sample performance measures in multi-server queueing systems can be
represented through formulas which are closely similar to those applied in
queueing networks. For instance, the average throughput may be defined exactly
as we have defined $  T_{K}^{n} $. To represent properly the remaining
measures, one however has to take into account the distinction between the
$k$th completion and the $k$th departure times, involved in the $  G/G/m  $
queue. With this distinction, replacing the symbols $  D_{k}  $ by
$  C_{k}  $ is required in the above formulas so as to provide appropriate
expressions for the sample performance measures of multi-server queues.

\subsection{Evaluation of System Performance}

We suppose now that the service times $  \tau_{k}  $ (and the interarrival
times $  \alpha_{k}  $ if they are given) are defined as random variables
$  \tau_{k} = \tau_{k}(\theta,\omega) $, where $  \theta \in \Theta  $ is a
set of decision parameters, and $  \omega  $ is a random vector. In this
case, as it results from the above representations of the queueing systems and
their performance, the arrival epochs $  A_{k}  $ and the departure epochs
$  D_{k} $, together with the sample performance measures also present random
variables. Let $  F_{K} = F_{K}(\theta,\omega) $, be a sample performance
measure of the system. As is customary, we define the system performance
measure associated with $  F_{K}  $ by the expected value
$$ \bm{F}_{K}(\theta) = E_{\omega}[F_{K}(\theta,\omega)]. $$

Based on a finite observation period, $  \bm{F}_{K}  $ is generally
referred to as {\em finite-horizon} performance measure. Another criterion,
{\em a steady-state} performance measure, intended to describe a long time
behaviour of a system is defined as
$$ \bm{F}(\theta)
    = \lim_{K \rightarrow \infty} E_{\omega}[F_{K}(\theta,\omega)]. $$

Although we may express sample performance measures in closed form, in the
case of general random variables determining the service times of customers,
it is often very difficult or even impossible to obtain analytically the
criteria $  \bm{F}_{K}  $ and especially $  \bm{F} $. In this
situation, one generally applies a simulation technique which allows of
obtaining values of $  F_{K}(\theta,\omega) $, and then estimates the system
performance by using the Monte Carlo approach. Note however, that information
concerning the explicit form of the sample performance measures normally
proves to be very useful to the simulation study and optimization of queueing
systems.

\section{Application of the Models}

In this section we briefly outline a selection of the application areas of the
recursive representation in simulation based analysis and optimization of
queueing systems. The section concludes with remarks concerning limitations on
the use of the models in representing queueing systems.

\subsection{Design of Simulation Algorithms}

Since recursive equations determine a global structure of changes in queueing
systems consecutively in a very natural way, they provide the basis for the
development of very efficient simulation procedures (see, e.g., \cite{Chen90,Greenberg91,Ermakov94}). Although the
simulation technique based on recursive representations of queueing systems
may rank below the traditional event-scheduling approach in its versatility,
the algorithms applying this technique are normally superior to others in
reducing time and memory costs. Moreover, these algorithms are usually best
suited to the implementation on parallel and vector processors. As an
illustration, one can consider parallel simulation algorithms in \cite{Greenberg91,Ermakov94}.

\subsection{Variance Reduction in Simulation}

Closely related to the queueing system simulation procedures are variance
reduction techniques which are intended to improve the accuracy of simulation
output (\cite{Cheng86}). In order for a variance reduction method to be
successfully employed in estimating a system performance
$  \bm{F}_{K} $, certain conditions normally have to be imposed on its
associated sample performance measure $  F_{K} $. Specifically, {\em the
antithetic variates\/} method and {\em the common random numbers\/} method
require that $  F_{K}  $ as a function of the random argument $  \omega  $
be monotone (see, e.g., \cite{Cheng86}). Examples of establishing such
monotonicity properties from the recursive representation of queueing systems
can be found in \cite{Rubinstein85}.

\subsection{Investigation of System Performance Measures}

Another area of applications of the models includes the investigation of
properties inherent in performance measures of queueing systems, such as
monotonicity and convexity with respect to system parameters $  \theta $. It
is normally not difficult to examine these properties for the systems
described by equations involving only the operations of maximum and addition
(e.g., tandem queues with both infinite and finite buffers). One can find an
extended discussion of this subject in \cite{Shanthikumar89a,Shanthikumar89b,Hu92}.

\subsection{Sensitivity Analysis and Estimation}

Since there are generally no explicit representations as functions of system
parameters $  \theta  $ available for the performance measure, one may
evaluate its sensitivity (or its gradient, when the parameters are continuous)
by no way other than through the use of estimates obtained from simulation
experiments. Very efficient procedures of obtaining gradient estimates may be
designed using new technique called {\em infinitesimal perturbation
analysis\/} (IPA) (see, e.g., \cite{Ho87}). The IPA algorithms which are actually
based on the recursive representations of queueing systems, can serve as an
important line of the application of the models under discussion (\cite{Cassandras88,Krivulin93}). Finally, these models provided a
useful framework for examining unbiasedness and consistency of IPA estimates
in \cite{Cassandras91,Hu92,Krivulin90,Krivulin93}.

\subsection{Limitations on the Use of the Models}

One can see that the general model of the network, as it has been presented
above, treats of queueing systems from the viewpoint of service facilities
rather than of particular customers. Specifically, for each node $  n  $
only the arrival and departure instants are essential, whereas it makes no
difference which of the customers proves to arrive or to depart. Moreover,
both times $  A_{k}^{n}  $ and $  D_{k}^{n}  $ do not need to be
associated with a single customer, as it normally happens in nodes with many
servers operating in parallel. As a consequence, the models do not allow of
representing systems with many classes of customers through recursive
equations in closed form. Finally, since nodes do not distinguish among
customers in some sense, the order in which customers are selected from a
queue for service is of no concern, and therefore, these models are incapable
of identifying distinct queue disciplines.

\bibliographystyle{utphys}

\bibliography{Recursive_equations_based_models_of_queueing_systems}

\end{document}